\documentclass[english]{article}
\usepackage[T1]{fontenc}
\usepackage{tipa}
\usepackage{tipx}
\usepackage{amsmath}
\usepackage{amsthm}
\usepackage{amssymb}
\usepackage{graphicx}
\usepackage[authoryear]{natbib}

\makeatletter

\providecommand{\tabularnewline}{\\}

\theoremstyle{definition}
\newtheorem{defn}{\protect\definitionname}
\theoremstyle{plain}
\newtheorem{thm}{\protect\theoremname}
\ifx\proof\undefined\
  \newenvironment{proof}[1][\proofname]{\par
    \normalfont\topsep6\p@\@plus6\p@\relax
    \trivlist
    \itemindent\parindent
    \item[\hskip\labelsep
          \scshape
      #1]\ignorespaces
  }{%
    \endtrivlist\@endpefalse
  }
  \providecommand{\proofname}{Proof}
\fi
\theoremstyle{plain}
\newtheorem{cor}{\protect\corollaryname}
\theoremstyle{definition}
\newtheorem{example}{\protect\examplename}

\usepackage{accents}

\makeatother

\usepackage{babel}
\providecommand{\corollaryname}{Corollary}
\providecommand{\definitionname}{Definition}
\providecommand{\examplename}{Example}
\providecommand{\theoremname}{Theorem}

\begin{document}
\title{Generalized Estimation and Information}
\author{Paul Vos\thanks{East Carolina University, vosp@ecu.edu} and Qiang
Wu\thanks{East Carolina University, wuq@ecu.edu}}

\maketitle
\global\long\def\tbar{\text{\ensuremath{\underbar{\ensuremath{\theta}}}}}%
\global\long\def\tbart{\text{\ensuremath{\underbar{\ensuremath{\theta}}}}^{\mathsf{t}}}%
\global\long\def\kbar{\underline{k}}%
\global\long\def\sbar{\underset{\bar{}}{s}}%
\global\long\def\partilde{\undertilde{\partial}}%
\global\long\def\ttilde{\undertilde{\theta}}%
\global\long\def\Tbar{\text{\ensuremath{\underbar{\ensuremath{\Theta}}}}}%
\global\long\def\nablatilde{\widetilde{\nabla}}%
\global\long\def\stilde{\undertilde{s}}%

\begin{abstract}
This paper extends the idea of a generalized estimator for a scalar
parameter \citep{Vos2022b} to multi-dimensional parameters both with
and without nuisance parameters. The title reflects the fact that
generalized estimators provide more than simply another method to
find point estimators, and that the methods to assess generalized
estimators differ from those for point estimators. By \emph{generalized
estimation} we mean the use of generalized estimators together with
an extended definition of \emph{information} to assess their inferential
properties. We show that Fisher information provides an upper bound
for the information utilized by an estimator and that the score attains
this bound. 
\end{abstract}

\section{Introduction}

 The maximum likelihood estimator need not be efficient and, among
the class of biased estimators, it need not be admissible. These issues
with maximum likelihood estimation and the parameter dependence of
other point estimators are addressed using generalized estimators.
Generalized estimators are described by information rather than variance
and the Fisher information provides an upper for the information of
an estimator. This bound applies to all generalized estimators; it
does not require estimators to be unbiased. The score is a generalized
estimator and its information equals the Fisher information. 

A point estimator assigns to each value $y$ in the sample space a
point in the parameter space $\Theta$. A generalized estimator $g$
assigns to each $y$ a function $g_{y}$ on $\Theta$ where $g_{y}(\theta)$
indicates the consistency of $\theta$ with $y$. The function $g_{y}$
can be thought of as a continuum of tests statistics evaluated at
$y$. The information of $g$ describes the average rate at which
these test statistics change with $\theta$. Section \ref{sec:One-Parameter-Family}
presents the scalar parameter case in a manner for natural extension
to multi-dimensional parameters in Section \ref{sec:Multi-parameter-Families}.
Section \ref{sec:Examples} presents two examples: one to illustrate
the role of information in assessing estimators and the other to illustrate
how confidence intervals can be obtained from a generalized estimate.

\section{One Parameter Families\label{sec:One-Parameter-Family}}

As we want inferences to be unaffected by the choice of parameter,
we describe the basics of inference without these. Parameterization
will be introduced to describe the smooth structures of estimators.

Let $M_{\mathcal{\mathcal{X}}}$ be a family of probability measures
having common support $\mathcal{X}$. While $\mathcal{X}$ can be
an abstract space, for most applications $\mathcal{X}\subset\mathbb{R}^{d}$.
Points in $M_{\mathcal{X}}$ serve as models for a population whose
individuals take values in $\mathcal{X}$.  We consider inference
for models from $M_{\mathcal{X}}$ based on a sample that is denoted
by $y$ and let $\mathcal{Y}$ be the corresponding sample space.
The relationship between $\mathcal{X}$ and $\mathcal{Y}$ will depend
on the sampling plan, conditioning, and dimension reduction using
sufficient statistics. For a simple random sample of size $n$ without
conditioning and no dimension reduction $\mathcal{Y}=\mathcal{X}^{n}$. 

Let $M=M_{\mathcal{Y}}$ be the family of probability measures obtained
from $M_{\mathcal{X}}$ using a sampling plan whose sample space is
$\mathcal{Y}$. For $\mathcal{Y}=\mathcal{X}^{n}$ 
\[
M=\left\{ m:m(y)=\prod m_{\mathcal{X}}(x_{i}),\ m_{\mathcal{X}}\in M_{\mathcal{X}}\right\} .
\]
 For the Bernoulli family of distributions, $\mathcal{X}=\left\{ 0,1\right\} $,
\[
M_{\mathcal{X}}=\left\{ m:0<m(1)<1,m(0)+m(1)=1\right\} .
\]
For a sample of size $n$ we use the sufficient statistic $y=\sum x_{i}$
so that $\mathcal{Y}=\left\{ 0,1,2,\ldots,n\right\} $ and
\begin{equation}
M=\left\{ m:m(y)={n \choose y}m_{\mathcal{X}}(1)^{y}m_{\mathcal{X}}(0)^{n-y},\ m_{\mathcal{X}}\in M_{\mathcal{X}}\right\} .\label{eq:MBin}
\end{equation}
When $\mathcal{Y}$ is open it will be convenient to let $m$ be a
probability density with respect to a dominating measure $\mu$. For
$\mathcal{X}=\mathbb{R}$ and function $\phi>0$ such that $\int\phi(x)d\mu=1$
there is a location family 
\[
M_{\mathcal{X}}=\left\{ m:m(x)=\phi(x-a),\ a\in\mathbb{R}\right\} .
\]
For a simple random sample with $y=\left(x_{1},x_{2},\ldots,x_{n}\right)^{{\tt t}}$
\[
M=\left\{ m:m(y)=\prod\phi(x_{i}-a),\ a\in\mathbb{R}\right\} .
\]
If $\phi(x)=\left(2\pi\right)^{-1/2}\exp\left(-\frac{1}{2}x^{2}\right)$
then $M_{\mathcal{X}}$ is the normal location family with unit variance.
Using the sufficient statistic $y=\bar{x}=(\sum x_{i})/n\in\mathcal{Y}=\mathbb{R}$,
\begin{equation}
M=\left\{ m:m\left(y\right)=\sqrt{n}\phi\left(\sqrt{n}\left(y-a\right)\right),a\in\mathbb{R}\right\} .\label{eq:normal}
\end{equation}
If $\phi(x)=\pi^{-1}\left(1+x^{2}\right)^{-1}$ then $M_{\mathcal{X}}$
is the Cauchy location family with unit scale factor. There is no
sufficient statistic of dimension less than $n$ so we use $y=\left(x_{1},x_{2},\ldots,x_{n}\right)^{{\tt t}}$,
\begin{equation}
M=\left\{ m:m(y)=\pi^{-n}\prod\left(1+\left(x_{i}-a\right)^{2}\right)^{-1},a\in\mathbb{R}\right\} .\label{eq:Cauchy}
\end{equation}

For real-valued measurable function $h$ we define the expected value
of $h$ at $m$,
\[
E_{m}h=\int_{\mathcal{Y}}h(y)m(y)d\mu
\]
when $\mathcal{Y}$ is open and $E_{m}h=\sum_{y\in\mathcal{Y}}h(y)m(y)$
when $\mathcal{Y}$ is discrete. We use the following Hilbert space
\[
H_{M}=\left\{ h:E_{m}h^{2}<\infty,\ \forall\ m\in M\right\} 
\]
which has a family of inner products indexed by $M$, 
\[
\langle h,h'\rangle_{m}=E_{m}\left(hh'\right)\mbox{for all }h,h'\in H_{M}.
\]
 When $E_{m}(hh')=0$ the vectors $h$ and $h'$ are \emph{$m$-orthogonal}
and we write $h\perp_{m}h'$. At each $m\in M$ there is a copy of
$H_{M}$ and this collection we denote by 
\[
H\!M=M\times H_{M}.
\]
The copy of $H_{M}$ at $m$ with inner product $\langle\cdot,\cdot\rangle_{m}$
is $H_{m}$ which we also write as $H_{m}M$ to indicate its relationship
to $H\!M$. For inference, $H_{m}$ will be restricted to the orthogonal
complement of the constant functions, $H_{m}^{\perp}=\{h\in H_{m}:E_{m}h=0\}$,
so that
\begin{equation}
H_{m}=H_{m}^{\perp}\oplus H_{m}^{0}\ \mbox{and}\ H_{m}^{\perp}\perp_{m}H_{m}^{0}.\label{eq:HHperpH0}
\end{equation}
Note $E_{m}h=\langle h,1\rangle_{m}$ and $H_{m}^{0}$ does not depend
on $m$. Since (\ref{eq:HHperpH0}) holds for each $m$ we write
\begin{equation}
H\!M=H^{\perp}\!M\oplus H^{0}M\ \mbox{and}\ H^{\perp}M\perp H^{0}M\label{eq:HHperpH02}
\end{equation}
where $\perp$ indicates $\perp_{m}$ holds for $H_{m}^{\perp}M=H_{m}^{\perp}$. 

As the notation suggests, $H\negmedspace M$ is a vector bundle on
$M$ with vector space $H_{M}$. It extends the tangent bundle $T\!M$
since $T\!M\subset H^{\perp}\!M$.  

For inference regarding models in $M$, we consider functions $g_{M}:\mathcal{Y}\times M\rightarrow\mathbb{R}$
such that 
\begin{equation}
g_{M}\left(\cdot,m\right)\in H_{m}^{\perp}\mbox{ for all }m\in M.\label{eq:hHbar}
\end{equation}
We also want $g_{M}$ to be a continuous on $M$, 
\begin{equation}
g_{M}\left(y,\cdot\right)\in C(M)\mbox{ for a.e. }y\in\mathcal{Y},\label{eq:hM}
\end{equation}
so that the expectation of $g_{M}$ is a continuous function. For
point estimators of a parameter, say $\theta$, the expectation of
the estimator $\hat{\theta}$ is a real number. To emphasize this
distinction we use the sans serif font to indicate the expectation
of $g_{M}$
\[
\mathsf{E}g_{M}\in C(M)\mbox{ while }E\hat{\theta}\in\mathbb{R}.
\]
Expectation $\mathsf{E}$ operates on $C(M)$-valued distributions,
whereas $E$ operates on $\mathbb{R}$-valued distributions. To be
a generalized estimator, $g_{M}(y,\cdot)$ will be required to have
continuous derivatives on $M$ and these will be described using parameterizations
that are diffeomorphisms.

We assume $M$ is a 1-dimensional smooth manifold. While more general
manifolds can be considered (e.g., Fisher's circle model), we will
only consider families that have a global parameterization 
\begin{equation}
\theta:M\rightarrow\Theta\subset\mathbb{R}\label{eq:thetaMap1dim}
\end{equation}
and are connected so that $\Theta$ is an open interval. For $g_{M}:\mathcal{Y}\times M\rightarrow\mathbb{R}$
we define $g_{\Theta}=g_{M}\circ\theta^{-1}:\mathcal{Y}\times\Theta\rightarrow\mathbb{R}$.
Unless more than one parameterization is being used, we drop the subscript
and write $g$ for $g_{\Theta}$. The\emph{ log likelihood function
on $\Theta$ for $y$} is the function defined by 
\[
\ell=\ell(y,\cdot)=\ell_{M}(y,\cdot)\circ\theta^{-1}
\]
where $\ell_{M}\left(y,m\right)=\log m\left(y\right)$. The \emph{score
function on $\Theta$ for $y$ }is
\[
s=\nabla\ell=\partial\ell/\partial\theta.
\]
We only consider $M$ such that $s(\cdot,\theta)\in H_{M}$ for all
$\theta\in\Theta$. Because $M$ is a smooth manifold $s\left(y,\cdot\right)\in C^{1}(\Theta)\ a.e.\ y$
and since $\mathsf{E}s=0$,
\begin{equation}
s(\cdot,\theta)\in H_{\theta}^{\perp}.\label{eq:Hm-1}
\end{equation}

These properties of $s$ are used to define generalized estimators. 
\begin{defn}
\label{def:generalized_estimator_1dim}A \emph{generalized estimator}
for scalar parameter $\theta$ is a function 
\[
g:\mathcal{Y}\times\Theta\longrightarrow\mathbb{R}
\]
and $g=g(y,\cdot)$ is the corresponding \emph{generalized estimate}
at $y$ if
\begin{align*}
\mbox{(i)} & \ \ g\left(y,\cdot\right)\in C^{1}(\Theta)\mbox{\ a.e.}\ y\\
\mbox{(ii)} & \ \ g\left(\cdot,\theta\right)\in H_{\theta}^{\perp}\mbox{ for all }\theta\\
\mbox{(iii)} & \ \ \mathsf{V}\left(g\right)>0
\end{align*}
where $\mathsf{V}\left(g\right)=\mathsf{E}\left(g^{2}\right)\in C^{1}(\Theta).$ 
\end{defn}
The space of generalized estimators for $\theta$ is $\mathcal{G}$
which we write as $\mathcal{G}_{\Theta}$ if we consider more than
one parameterization. Any function $f\not\in H_{\theta}^{\perp}$
that satisfies $f\in H_{\theta}$ and conditions (i) and (iii) of
Definition \ref{def:generalized_estimator_1dim} is a \emph{pre generalized
estimator}, or simply, a \emph{pre estimator}. For any pre-estimator
$f$, its \emph{orthogonalization }
\begin{equation}
f^{\perp}=f-f^{\top}\in\mathcal{G}.\label{eq:preEst1dim}
\end{equation}
where $f^{\top}=\mathsf{E}f$. 

\citet{Godambe1960} has similar criteria but allows $\mathsf{V}(g)=0$
for some $\theta\in\Theta$ and adds that $\mathsf{E}\left(\nabla g\right)^{2}>0$
so that $\mathsf{E}\left(\nabla g\right)$ can never be zero on $\Theta$.
We do not need this restriction since we describe estimators in terms
of information rather than variance. Allowing $\mathsf{E}\left(\nabla g\right)$
to be zero will be useful for nuisance parameters in the multi-dimension
setting. Because $\mathsf{V}(g)>0$ we can define the\emph{ standardization
of} $g$ as
\[
\bar{g}=\frac{g}{\sqrt{\mathsf{V}(g)}}.
\]
Since $\bar{g}(\cdot,\theta)\in H_{\theta}^{\perp}$ is a vector of
unit length, $\bar{g}$ is also called the \emph{direction} of $g$.
Standardized estimators are the same in every parameterization. That
is, for any $m'\in M$, $\bar{g}_{\Theta}(\cdot,\theta')=\bar{g}_{\Xi}(\cdot,\xi')$
where $\theta'=\theta(m')$ and $\xi'=\xi(m')$. 

A non-degenerate point estimator $\hat{\theta}$ whose first two moments
are smooth functions on $\Theta$ is a pre-estimator so that 
\[
\hat{\theta}-\mathsf{E}\hat{\theta}\in\mathcal{G}.
\]
We use the sans serif notation because as a pre-estimator $\hat{\theta}(y,\cdot)$
is a function on the parameter space, the constant function taking
the value of the point estimate at $y$. The estimator need not be
unbiased, so that generalized estimation can be used to compare biased
and unbiased point estimators as well as estimators not constrained
to be constant on $\Theta$. Generalized estimators are compared in
terms of their information. 
\begin{defn}
\label{def:Lambda}The \emph{information for scalar parameter $\theta$
utilized by }$g$ is 
\begin{equation}
\Lambda(g)=\left(\mathsf{E}\nabla\bar{g}\right)^{2}=\frac{(\mathsf{E}\nabla g)^{2}}{\mathsf{E}(g^{2})}.\label{eq:squaredSlopeDefined}
\end{equation}
where the second equality follows from the definition of $\bar{g}$
and $\mathsf{E}g=0$.
\end{defn}
The Fisher information for a sample of size $n$, $I_{(n)}$, and
the Fisher information in a single observation, $I_{(1)}$, satisfy
$I_{(n)}=nI_{(1)}$. This relationship also holds for the information
utilized by an estimator
\begin{equation}
\Lambda(g_{(n)})=n\Lambda(g_{(1)}).\label{eq:Lambda(n) Lambda(1)}
\end{equation}
When considering only samples of size $n$ we use $I=I_{(n)}$ and
$\Lambda(g)=\Lambda(g_{(n)})$. 

As the score is the archetype for a generalized estimator $g$, the
log likelihood function is the archetype for the scalar potential. 
\begin{defn}
\label{def:scalar potential}A \emph{scalar potential} \emph{of} $g$
is any function $G:\mathcal{Y}\times\Theta\longrightarrow\mathbb{R}$
such that $\nabla G=g$.
\end{defn}
While $G\not\in\mathcal{G}$ we define the information utilized by
$G$ to be the information of its derivative: $\Lambda(G)=\Lambda(g)$.
Information is a local property and so does not distinguish between
a generalized estimator and its scalar potential. The scalar potential
is useful for finding confidence regions especially when the parameterization
is multidimensional.

We assume differentiation commutes with the integral sign so for any
pre-estimator $f$
\begin{equation}
\nabla\left(\mathsf{E}f\right)=\mathsf{E}\left(\nabla f\right)+\left(\nabla\mathsf{E}\right)\left(f\right)\label{eq:NablaEcommute}
\end{equation}
where $\left(\nabla\mathsf{E}\right)$ is the linear operator on $H_{M}$
defined by 
\[
\left(\nabla\mathsf{E}\right)(h)=\mathsf{E}\left(\left(\nabla\ell\right)h\right).
\]
Note that we use $f$ and $g$ for functions on $\mathcal{Y}\times\Theta$
while $h\in H_{M}$ is a function on $\mathcal{Y}$. For generalized
estimator $g$, $\mathsf{E}g$ vanishes so (\ref{eq:NablaEcommute})
becomes, after switching left- and right-hand sides, the \emph{score
equation} 
\begin{equation}
\mathsf{E}\left(\nabla g\right)+\mathsf{E}\left(sg\right)=0.\label{eq:scoreIdentity1}
\end{equation}
When $g=s$, the score equation gives the equivalent definitions of
the Fisher information for $\theta$ 
\[
I=-\mathsf{E}(\nabla s)=\mathsf{E}(s^{2}).
\]

The information upper bound follows from the score identity. 
\begin{thm}
\label{thm:1}The information for $\theta$ utilized by $g$ is bounded
by the Fisher information 
\begin{eqnarray*}
\Lambda\left(g\right) & \le & I.
\end{eqnarray*}
Furthermore, the score $s$ attains this bound and for any $g\in\mathcal{G}$
\begin{eqnarray*}
\Lambda(g) & = & \mathsf{V}(\mathsf{P}_{g}s)\\
 & = & \mathsf{R}^{2}I
\end{eqnarray*}
where $\mathsf{P}_{g}s$ is the projection of $s$ onto the space
spanned by $g$ and $\mathsf{R}=\mathsf{E}(\bar{s}\bar{g})$ is the
correlation between $s$ and $g$. 
\end{thm}
\begin{proof}
From the score equation
\begin{equation}
\Lambda(g)=\mathsf{E}^{2}(s\bar{g}).\label{eq:proof1}
\end{equation}
The second displayed equality follows upon noting that $\mathsf{E}^{2}(s\bar{g})=\mathsf{E}^{2}(\bar{s}\bar{g})I=\mathsf{R}^{2}I$.
The first equality follows by expressing the projection using basis
vector $\bar{g}$
\begin{eqnarray*}
\mathsf{V}(\mathsf{P}_{g}s) & = & \mathsf{V}\left(\mathsf{E}(s\bar{g})\bar{g}\right)\\
 & = & \mathsf{E}^{2}(s\bar{g}).
\end{eqnarray*}
\end{proof}
Efficiency of a point estimator is defined using the ratio of its
variance to the variance bound. Efficiency of a generalized estimator
is defined as the ratio of its information to the information bound,
$I$. 
\begin{defn}
The \emph{$\Lambda$-efficiency of $g$ is }
\[
\mbox{Eff}^{\Lambda}(g)=I^{-1}\Lambda(g).
\]
\end{defn}
An immediate corollary of Theorem \ref{thm:1} is that the $\Lambda$-efficiency
is the square of the correlation between the estimator and the score. 
\begin{cor}
\label{cor:1}

\begin{eqnarray*}
\mathrm{Eff}^{\Lambda}(g) & = & \mathsf{V}\left(\mathsf{P}_{g}\bar{s}\right)\\
 & = & \mathsf{R}^{2}.
\end{eqnarray*}
\end{cor}
The $\Lambda$-efficiency of a point estimator $\hat{\theta}$ is
the $\Lambda$-efficiency of its generalized estimator $g_{\hat{\theta}}=\hat{\theta}-\mathsf{E}\hat{\theta}$.
When $\hat{\theta}$ is unbiased $\Lambda(g_{\hat{\theta}})=\mathsf{V}^{-1}(\hat{\theta})$
so that $\Lambda$-efficiency is identical to efficiency based on
variance. 

Even though these efficiencies can take the same numerical value,
it is incorrect to characterize the information as the reciprocal
of the variance. The information at $\theta'=\theta(m')$, $\Lambda(g)|_{\theta=\theta'}$,
is a measure of how $g$ changes in a neighborhood $m'\in M$; that
is, information depends on $M$. The variance at $\theta'$, $\mathsf{V}(g)|_{\theta=\theta'}$,
depends only on $m'$; it is the same for the countless manifolds
we could choose that contain $m'$.  Another difference is that variance
is defined on horizontal distributions while information is defined
on vertical distributions. Horizontal and vertical distributions are
described in Example \ref{exa:binom20}. 

\begin{example}
\label{exa:binom20}We consider inference for the proportion of a
population having a genetic variation or other specified characteristic.
We let $1$ (0) indicate the characteristic is present (absent) so
$\mathcal{X}=\left\{ 0,1\right\} $ and for a sample of size $n$,
$M$ is given by (\ref{eq:MBin}).  Figure \ref{Figure-1: Bin20 score}
shows the standardized score
\[
\bar{s}=\frac{y-np}{\sqrt{np(1-p)}}
\]
where $n=20$ and $p$ is the parameter defined by $p(m)=m(1)$ with
parameter space $P=(0,1)$. The graph of the estimate $\bar{s}_{y}$
when $y=6$ is the black curve. The estimator $\bar{s}$ is represented
by the family of 21 curves, one for each $y$ in the sample space
(unrealized estimates are shown in white). 

\begin{figure}
\includegraphics[width=5in,height=5in,keepaspectratio]{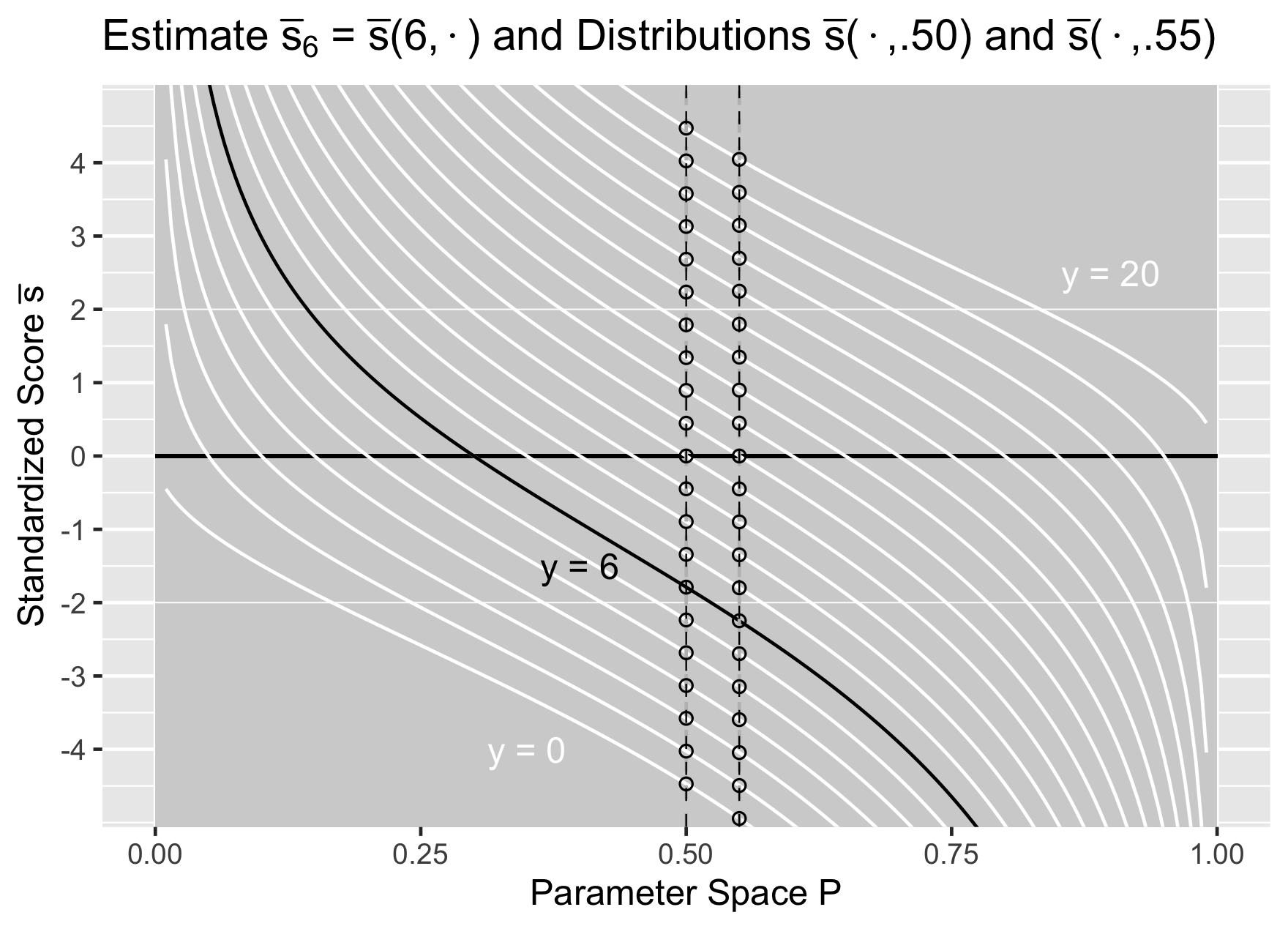}\caption{\label{Figure-1: Bin20 score}The standardized score estimate $\bar{s}_{6}$
obtained from the sample with $y=6$ and $n=20$ for the Bernoulli
manifold with the parameter $p=m(1)$ is shown by the black curve.
The standardized score estimator $\bar{s}$ is represented by the
family of 21 curves, one for each $y$ in the sample space (unrealized
estimates are shown in white). Of the continuum of vertical slices
two are shown at $p=.50$ and $p=.55$. The distribution of the point
estimate $\hat{p}$ is shown by the intersection of these 21 curves
with the horizontal axis. Note that for two of these curves the intersection
occurs for a value outside of the parameter space. }

\end{figure}

Of the continuum of vertical slices two are shown, one at $p=.50$
and another at $p=.55$. Every vertical slice for $0<p<1$ intersects
all 21 curves and while the ordinate of these points of intersection
depends on $p$ the resulting distributions all have mean zero and
variance one. These vertical distributions are the same in every parameterization.
For any parameter $\theta$, $\bar{s}(y,p(m'))=\bar{s}_{\Theta}(y,\theta(m'))$
for all $y$ and all $m'$. In contrast, the abscissa values obtained
from the intersection of these curves with the parameter axis are
the same for all $p$ but the mean and variance of these horizontal
distributions depends on the value of the parameter and on the choice
of parameterization. The horizontal distributions describe the inferential
properties in terms of the mean and variance of the roots of $s$
while the vertical distributions describe how each estimate $\bar{s}_{y}$
changes with the parameter. 

When the maximum likelihood estimator exists and is unique it is,
by definition, the parameter-intercept of the score, $\hat{p}=s^{-1}(0)$.
For $y=0$ and for $y=20$, the maximum likelihood estimate does
not exist since $s_{y}$ does not cross the parameter axis. Even when
the point estimate does not exist, confidence regions can be constructed
from the standardized score $\bar{s}$. All 21 estimates $\bar{s}_{y}$
provide $z$-standard deviation intervals
\[
\mbox{CI}_{-z}(y)=\left\{ p:\bar{s}_{y}(p)\ge-z\right\} ,\mbox{CI}_{+z}(y)=\left\{ p:\bar{s}_{y}(p)\le z\right\} .
\]
The intersection of the curve $\bar{s}_{6}$ with the white lines
at $\bar{s}_{y}=\pm2$ in Figure \ref{Figure-1: Bin20 score} show
the endpoints of $\mbox{CI}_{-2}(6)$ and $\mbox{CI}_{+2}(6)$. Since
generalized estimators are parameter invariant, these intervals correspond
to subsets of the space of models $M$. The interpretation of these
intervals can be stated in terms of their complement: if the true
model is not in $\mbox{CI}_{-2}(y)$ or $\mbox{CI}_{+2}(y)$ then
the score test for the observed data $y$ is at least two standard
deviations from zero. That is, for models outside these intervals
the observed data $y$ would be improbable since the score is at least
two standard deviations from zero. Intervals based on tail probabilities
can be obtained by allowing $z$ to be a function of the parameter;
for $\mbox{CI}_{+z}(6)$ the value for $z$ would be obtained using
the mass assigned to the values $\{0,1,\ldots,5,6\}$.

Figure \ref{fig:Twice-the-log} shows the log likelihood ratio statistic
$S$ for $y=6$ and its distribution on the other 20 values in the
sample space. The vertical slices at $p=.50$ and $p=.55$ correspond
to those from Figure \ref{Figure-1: Bin20 score} but the circles
are only plotted when the slope of the intersecting curve is negative.
Each vertical slice has 6 points of intersection corresponding to
samples as extreme as $y=6$. The resulting p-value is the same as
for the score. This will be true for any vertical slice so that inference
from the score and the signed log likelihood ratio are identical in
this example. This will not be true when the curves of the estimator
$g$ intersect. Also, inference from $g$ and unsigned scalar potential
function $G$ will not be identical. In particular, the score and
unsigned log likelihood ratio are not identical in this example.

\begin{figure}
\includegraphics[width=5in,height=5in,keepaspectratio]{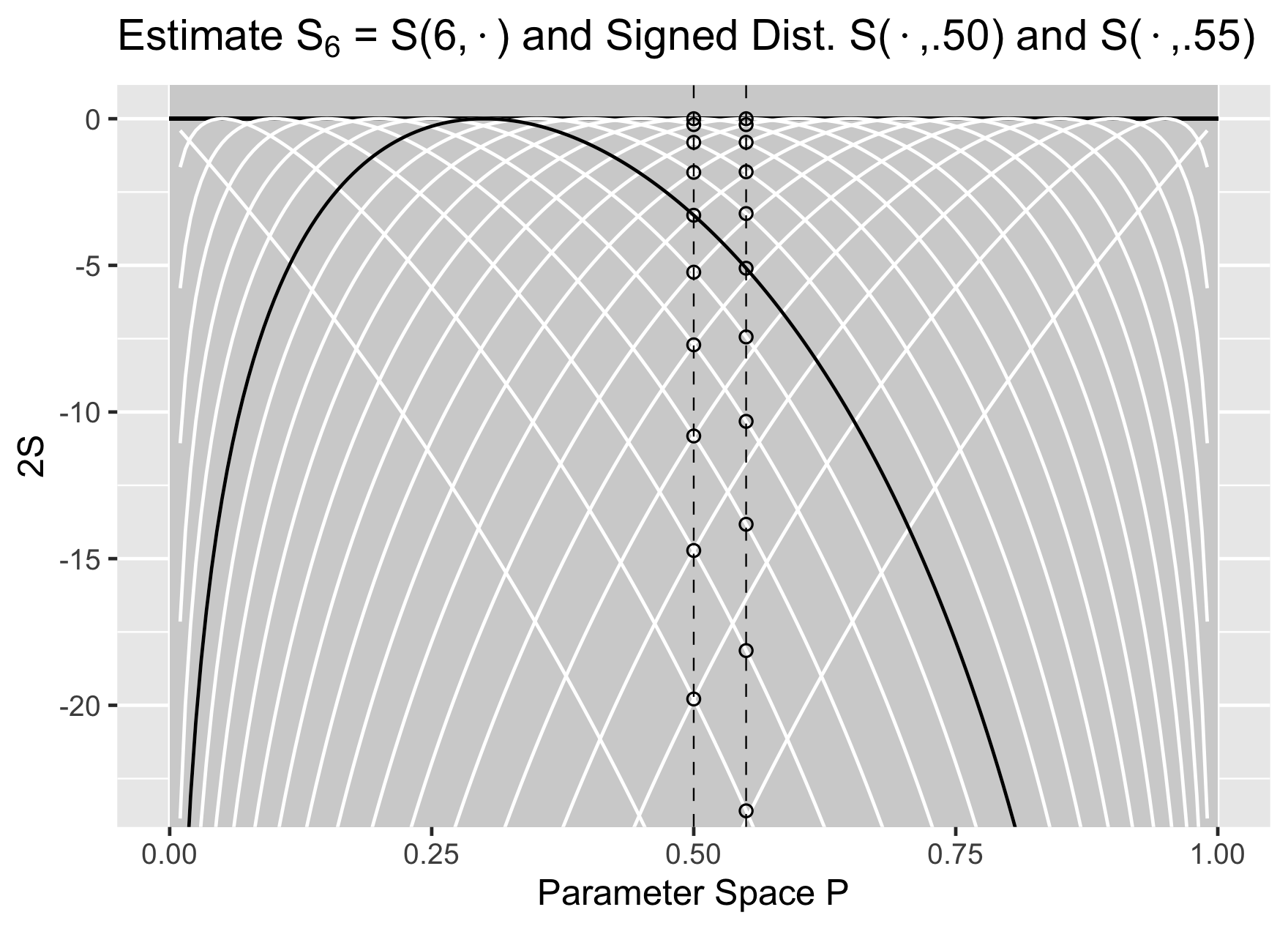}\caption{\label{fig:Twice-the-log}Twice the log likelihood ratio statistic
obtained from observing $y=6$ out of a sample of size $n=20$ for
the Bernoulli manifold with the parameter $p=m(1)$ is shown by the
black curve. The distribution of twice the log likelihood ratio statistic
is represented by the black curve and 20 white curves. }
\end{figure}
\end{example}
 
\begin{example}
\label{Cauchy}-- We consider the same population as before but
now the variable of interest is a measured quantity and we choose
$M_{\mathcal{X}}$ to be the Cauchy family so that for a random sample
of size $n$, $M$ is given by (\ref{eq:Cauchy}). For comparison
we also consider models from the Normal family for which the family
of sampling distributions is given by (\ref{eq:normal}); we use $M_{{\tt Gaus}}$
to identify this manifold. For parameterization $\theta$, the graph
of a generalized estimate $\bar{g}_{y}$ for an observation $y=\left(x_{1},x_{2},\ldots,x_{n}\right)^{{\tt t}}$
is a curve over the parameter space $\Theta$. This corresponds to
the black curve in the previous example. The distribution of the estimator
$\bar{g}$ is more difficult to represent since there are a continuum
of curves indexed by $y$. For $M_{{\tt Gaus}}$ there is also a continuum
of curves but now the sufficient statistic $\bar{x}=n^{-1}\sum x_{i}$
provides a one dimensional index. Nevertheless, the properties of
the vertical distributions for $M$ and $M_{{\tt Gaus}}$ still hold
and confidence regions for $g_{y}$ are defined in the same way. 

\end{example}

\section{Multi-parameter Families\label{sec:Multi-parameter-Families}}

We consider inference for a parameter $\theta=\left(\theta^{1},\theta^{2},\ldots,\theta^{k}\right)^{{\tt t}}\in\mathbb{R}^{k}$
in the presence of a $k'$-dimensional nuisance parameter $\undertilde{\theta}=(\undertilde{\theta}^{1},\undertilde{\theta}^{2},\ldots,\undertilde{\theta}^{k'})^{{\tt t}}$
so that $M$ is a manifold of dimension $(k+k')$ and $\tbar^{{\tt t}}=(\theta^{{\tt t}},\undertilde{\theta}^{{\tt t}})$
is a global parameterization $\tbar:M\rightarrow\Tbar$. We use $\overline{\nabla}$,
$\nabla$, and $\nablatilde$ to indicate differentiation with respect
to $\underline{\theta}$, $\theta$, and $\undertilde{\theta}$, respectively,
so that 
\[
\undertilde{s}=\nablatilde\ell=(\partial\ell/\partial\undertilde{\theta}^{1},\partial\ell/\partial\undertilde{\theta}^{2},\ldots,\partial\ell/\partial\undertilde{\theta}^{k'})^{{\tt t}}.
\]
Note that subscripts are used for the components of $g$ while superscripts
are used for $\theta$. This convention allows us to use the Einstein
summation convention for calculations involving bases. It also reminds
us that the component $g_{a}$ is not a point estimate for $\theta^{a}$;
if it were, we would want to use superscripts for the components of
$g$. While $\theta$ and $g$ are both $k$-tuples, geometrically
$\theta$ is a contra-variant (tangent) vector while $g$ is a covariant
vector as its components co-vary with the change of basis. 

Generalized estimators may depend on the value of the nuisance parameter
but we can make them independent of the nuisance \emph{parameterization}
by restricting to functions that are orthogonal to $\stilde$. For
any fixed $m_{\circ}\in M$ there is a $k'$-dimensional submanifold
through $m_{\circ}$
\[
M|_{m_{\circ}}=\left\{ m\in M:\theta(m)=\theta_{\circ}\right\} 
\]
where $\theta_{\circ}=\theta(m_{\circ}).$The tangent space of $M|_{m_{\circ}}$
at $m\in M|_{m_{\circ}}$ is
\[
\widetilde{T}_{m}M=\mbox{span}\{\stilde(\cdot,\tbar)\}|_{\underline{\theta}^{{\tt t}}=(\theta_{\circ}^{{\tt t}},\undertilde{\theta}^{{\tt t}})}.
\]
We will require estimators to be orthogonal to $\widetilde{T}_{m}M$
and so define
\[
H_{m}^{\bot}=\left\{ h\in H_{M}:E_{m}h=0,h\perp_{m}\widetilde{T}_{m}M\right\} .
\]

Equations (\ref{eq:HHperpH0}) and (\ref{eq:HHperpH02}) for the one
dimensional case become
\begin{eqnarray}
H_{m} & = & H_{m}^{\perp}\oplus\widetilde{T}_{m}M\oplus H_{m}^{0}\label{eq:HH_multivariate1}\\
H\!M & = & H^{\perp\!}M\oplus\widetilde{T}\!M\oplus H^{0}\!M\nonumber 
\end{eqnarray}
When $M$ is parameterized by $\text{\ensuremath{\underbar{\ensuremath{\theta}}}}$,
$m$ in (\ref{eq:HH_multivariate1}) is replaced with $\text{\ensuremath{\underbar{\ensuremath{\theta}}}=\ensuremath{\underbar{\ensuremath{\theta}}}}(m)$.

\begin{defn}
\label{def:genEstk}A \emph{generalized estimator} \emph{for a $k$-dimensional
parameter} $\theta$ is a function 
\[
g:\mathcal{Y}\times\underline{\Theta}\longrightarrow\mathbb{R}^{k}
\]
and $g_{y}=g(y,\cdot)$ is the corresponding \emph{generalized estimate}
at $y$ if
\begin{align*}
\mathrm{(I)} & \ \ g\left(y,\cdot\right)\in C^{1}(\underline{\Theta},\mathbb{R}^{k})\mbox{\ a.e.}\ y\\
\mathrm{(II)} & \ \ g(\cdot,\tbar)\in H_{\tbar}^{\bot}\mbox{ for all }\tbar\in\underline{\Theta}\\
\mbox{(III)} & \ \ \mathsf{V}\left(g\right)>0
\end{align*}
where $\mathsf{V}(g)=\mathsf{E}(gg^{\mathsf{t}})\in C^{1}\left(\underline{\Theta},\mathbb{R}^{k\times k}\right)$. 
\end{defn}
\vphantom{}

The space of generalized estimators for $\theta$ is $\mathcal{G}$
which we write as $\mathcal{G}_{\Theta}$ if we consider more than
one parameterization.  If $f_{\text{\ensuremath{\underbar{\ensuremath{\theta}}}}}=f(\cdot,\tbar)\in H_{\underline{\theta}}$
for all $\tbar\in\Tbar$ and satisfies conditions (I) and (III) of
Definition \ref{def:genEstk} but $f_{\tbar}\not\in H_{\underline{\theta}}^{\perp}$\emph{,
}then $f$ is a\emph{ pre-estimator}. The \emph{orthogonalization
of} $f$ \emph{at} $\text{\ensuremath{\underbar{\ensuremath{\theta}}}}$
\begin{equation}
f_{\text{\ensuremath{\underbar{\ensuremath{\theta}}}}}^{\bot}=f_{\text{\ensuremath{\underbar{\ensuremath{\theta}}}}}-f_{\underline{\theta}}^{\top}\in H_{\text{\ensuremath{\underbar{\ensuremath{\theta}}}}}^{\perp}\label{eq:fperp}
\end{equation}
where $f_{\underline{\theta}}^{\top}=E_{\underline{\theta}}(f_{\text{\ensuremath{\underbar{\ensuremath{\theta}}}}})+\widetilde{P}_{\text{\ensuremath{\underbar{\ensuremath{\theta}}}}}f_{\text{\ensuremath{\underbar{\ensuremath{\theta}}}}}$
and $\widetilde{P}_{\text{\ensuremath{\underbar{\ensuremath{\theta}}}}}$
is the orthogonal projection onto $\widetilde{T}_{\text{\ensuremath{\underbar{\ensuremath{\theta}}}}}\!M$.
Since (\ref{eq:fperp}) holds for all $\text{\ensuremath{\underbar{\ensuremath{\theta}}}}$
and expectation and orthogonal projections are smooth functions we
have
\[
f^{\perp}=f-f^{\top}\in\mathcal{G}
\]
where $f^{\top}=\mathsf{E}f+\widetilde{\mathsf{P}}f$.

The score $\nabla\ell$ is a pre-estimator so that we define $s$
to be the\emph{ orthogonalized score }
\[
s=(\nabla\ell)^{\bot}\in\mathcal{G}.
\]
The Fisher information for $\theta$ is $I=\mathsf{V}\left(\nabla\ell\right)$
and the nuisance orthogonalized Fisher information is $I^{\perp}=\mathsf{V}\left((\nabla\ell)^{\bot}\right)=\mathsf{V}\left(s\right)$;
both can be functions of $\undertilde{\theta}$ but only $I^{\perp}$
is the same for all nuisance parameterizations. 

The relationship between the score (information) and the orthogonalized
score (orthogonalized information) expressed in the $\text{\ensuremath{\underbar{\ensuremath{\theta}}}}$
parameterization is
\begin{eqnarray*}
\left(\nabla\ell\right)^{\perp} & = & \nabla\ell-I_{\nabla\widetilde{\nabla}}\undertilde{I}^{-1}\widetilde{\nabla}\ell\\
I^{\perp} & = & I-I_{\nabla\widetilde{\nabla}}\undertilde{I}^{-1}I_{\widetilde{\nabla}\nabla}
\end{eqnarray*}
where 
\[
\underline{I}=I_{\overline{\nabla}\overline{\nabla}}=\begin{pmatrix}\begin{array}{ll}
I & I_{\nabla\widetilde{\nabla}}\\
I_{\widetilde{\nabla}\nabla} & \undertilde{I}
\end{array}\end{pmatrix}
\]
and $I=I_{\nabla\nabla}$ and $\undertilde{I}=I_{\widetilde{\nabla}\widetilde{\nabla}}$
are the Fisher informations for $\theta$ and $\undertilde{\theta}$.
When $I_{\nabla\widetilde{\nabla}}$ vanishes on $\text{\ensuremath{\underbar{\ensuremath{\Theta}}}}$,
parameterizations $\theta$ and $\undertilde{\theta}$ are \emph{orthogonal}.

The definition of the scalar potential in the multi-parameter case
is straight forward. And, as in the scalar parameter case, the log
likelihood $\ell$ is the scalar potential for $s$. 
\begin{defn}
\label{def:scalar potential-1}A \emph{scalar potential} \emph{of}
$g\in\mathcal{G}$ is any function $G:\mathcal{Y}\times\text{\ensuremath{\underbar{\ensuremath{\Theta}}}}\longrightarrow\mathbb{R}$
such that $g=(\nabla G)^{\bot}$.
\end{defn}
The multivariate version of (\ref{eq:NablaEcommute}) is
\begin{equation}
\nabla\mathsf{E}(f^{{\tt t}})=\mathsf{E}\left(\nabla f^{{\tt t}}\right)+\left(\nabla\mathsf{E}\right)(f^{{\tt t}})\label{eq:NablaEcommuteMult}
\end{equation}
where 
\[
\left(\nabla\mathsf{E}\right)\left(f^{{\tt t}}\right)=\mathsf{E}\left(\left(\nabla\ell\right)f^{{\tt t}}\right).
\]
Since $g\in H_{M}^{\perp}$ we have $\mathsf{E}\left(\left(\nabla\ell\right)g^{{\tt t}}\right)=\mathsf{E}\left(sg^{{\tt t}}\right)$
so that the multivariate version of the \emph{score equation} (\ref{eq:scoreIdentity1})
is
\begin{equation}
\mathsf{E}\left(\nabla g^{{\tt t}}\right)+\mathsf{E}\left(sg^{{\tt t}}\right)=0.\label{eq:scoreIdentk}
\end{equation}
Differentiating with respect to the nuisance parameter we obtain
\begin{equation}
\mathsf{E}(\widetilde{\nabla}g^{{\tt t}})+\mathsf{E}(\stilde g^{{\tt t}})=0\label{eq:Enablanuis=00003D0}
\end{equation}
so that $g$ being nuisance orthogonal means that the average change
of $g$ in the direction of the nuisance parameter is zero. 

For the mean slope to be meaningful we need to use its standardized
version. 
\begin{defn}
\label{def:gbark} For $g\in\mathcal{G}$, define 
\[
\bar{g}=\mathsf{V}^{-1/2}g
\]
where $\mathsf{V}=\mathsf{V}\left(g\right)$ so that $\mathsf{V}\left(\bar{g}\right)$
is $I_{\mathsf{id}}$, the $k\times k$ identity matrix. Any $g$
such that $\mathsf{V}\left(g\right)=I_{\mathsf{id}}$ is called a
\emph{standardized estimator}.
\end{defn}
\smallskip{}

\begin{defn}
\label{Lambdak}\emph{ The information for $\theta$ utilized by $g$
is 
\begin{align*}
\Lambda\left(g\right) & =\left(\mathsf{E}\nabla\bar{g}^{\mathsf{t}}\right)\left(\mathsf{E}\nabla\bar{g}^{\mathsf{t}}\right)^{\mathsf{t}}\\
 & =\left(\mathsf{E}\nabla g^{\mathsf{t}}\right)\mathsf{V}^{-1}(g)\left(\mathsf{E}\nabla g^{\mathsf{t}}\right)^{\mathsf{t}}.
\end{align*}
The }scalar\emph{ information for $\theta$ utilized by $g$ is
\[
\lambda(g)=\mathop{tr}\Lambda(g).
\]
}
\end{defn}
\smallskip{}

Note $\Lambda(g)\in C^{1}(\text{\ensuremath{\underbar{\ensuremath{\Theta}}}},\mathbb{R}^{k}\times\mathbb{R}^{k})$.
Using the Frobenius norm for matrix $A$, $||A||=\sqrt{\mathop{tr}(A^{{\tt t}}A)}$,
we see that the scalar information is the square of the norm of $\mathsf{E}\nabla\bar{g}^{{\tt t}}$
\[
\lambda(g)=||\mathsf{E}\nabla\bar{g}^{{\tt t}}||^{2}.
\]
By replacing $\nabla$ with $\widetilde{\nabla}$ in Definition \ref{Lambdak}
we could define $\undertilde{\Lambda}(g)$, the information for $\undertilde{\theta}$.
However, equation (\ref{eq:Enablanuis=00003D0}) shows $\undertilde{\Lambda}(g)=0$
for all $g\in\mathcal{G}$. Restricting estimators to be orthogonal
to the space spanned by the nuisance parameters makes inferences independent
of the choice of the nuisance parameter but also means that estimators
for the parameter of interest have no information for the nuisance
parameter. 
\begin{thm}
\label{thm:Lambda}For k-dimensional parameter $\theta$ let $s=(\nabla\ell)^{\perp}$
and let $I^{\perp}=\mathsf{V}(s)$ be the orthogonalized Fisher information
for $\theta$. For any $g\in\mathcal{G}$, $\Lambda(g)\le I^{\perp}$
and $s$ attains this bound, $\Lambda(s)=I^{\perp}$. Furthermore,
\begin{eqnarray*}
\Lambda(g) & = & \mathsf{V}(\mathsf{P}_{g}s)=\mathsf{V}(\mathsf{P}_{g}\nabla\ell)\\
 & = & (I^{\perp})^{1/2}\mathsf{R}\mathsf{R}^{{\tt t}}(I^{\perp})^{1/2}
\end{eqnarray*}
 where $\mathsf{R}=\mathsf{E}(\bar{s}\bar{g}^{{\tt t}})$ is the correlation
matrix between $s$ and $g$. 
\end{thm}
\begin{proof}
The displayed equations in the Theorem are obtained from the score
equation (\ref{eq:scoreIdentk}) which gives
\[
\Lambda(g)=\mathsf{E}(s\bar{g}^{{\tt t}})\mathsf{E}(\bar{g}s^{{\tt t}}).
\]
The first equation follows from the definition of the projection and
its variance: $\mathsf{P}_{g}s=\mathsf{E}(s\bar{g}^{{\tt t}})\bar{g}$
so $\mathsf{V}(\mathsf{P}_{g}s)=\mathsf{E}(s\bar{g}^{{\tt t}})\mathsf{E}(\bar{g}s^{{\tt t}})$.
The second equation follows because $\nabla\ell=s+(\nabla\ell)^{\top}$
and $g$ is orthogonal to $\left(\nabla\ell\right)^{\top}$. The third
equation follows from $\mathsf{E}(s\bar{g}^{{\tt t}})=(I^{\bot})^{1/2}\mathsf{E}(\bar{s}\bar{g}^{{\tt t}})$
since $\mathsf{V}(s)=I^{\bot}$. The inequality $\Lambda(g)\le I^{\bot}$
follows because the squared length of a projection cannot be longer
than the original vector.
\end{proof}
When there are no nuisance parameters Theorem \ref{thm:Lambda} holds
with $I^{\bot}=I$ and $s=\nabla\ell$. 
\begin{defn}
\label{def:effk}The $\Lambda$-\emph{efficiency of} $g$ is
\[
\mbox{Eff}^{\Lambda}\left(g\right)=(I^{\perp}){}^{-1/2}\Lambda(g)(I^{\perp}){}^{-1/2}.
\]
\end{defn}
Corollary \ref{cor:2} follows immediately from Theorem \ref{thm:Lambda}.
\begin{cor}
\label{cor:2}
\begin{eqnarray*}
\mathrm{Eff}^{\Lambda}\left(g\right) & = & \mathsf{V}(\mathsf{P}_{g}\bar{s})\\
 & = & \mathsf{R}\mathsf{R}^{{\tt t}}.
\end{eqnarray*}
\end{cor}

\section{Examples\label{sec:Examples}}

\subsection{Normal and $t$-distributions\label{subsec:ExampleNormal-and--distributions}}

We consider two one-dimensional manifolds: the normal family and the
family of $t$ distributions with 3 degrees of freedom. Both are location
families so Fisher information is the same at each distribution in
the manifold. We compare three estimators: the sample mean, sample
median and the mle obtained from the $t_{3}$ distribution. We did
not include the score for the $t_{3}$ distribution since it is very
close to the corresponding mle. The sample mean is the mle for normal
data.

The sample mean attains the information bound for the normal family
and the $t_{3}$ score attains the information bound for the $t_{3}$
family. We use the information of these estimators to assess the cost
of model misspecification and explore the relationship between information
and the tails of the distribution. 

Figure \ref{fig:Comparing-Estimators:-Mean,} is based on 100,000
samples of size 10 from a normal distribution and another 100,000
samples of size 10 from the $t_{3}$ distribution. For the graph on
the left, 99 quantiles, from .005 to .995, obtained from the 100,000
sample means for the normal data were calculated. Using the empirical
cdf for the 100,000 medians these 99 quantiles gave 100 tail areas
(the median was included in both tails). Each tail area $T\!A$ was
converted to a $\zeta$-score that measures the distance into the
tail of the distribution. For continuous random variable $X$ define
$\zeta:X\rightarrow\mathbb{R}$ by 
\[
\zeta=\left\{ \begin{array}{ll}
\log_{2}(2\mbox{Pr}(X\le x)) & \mbox{if }\mbox{Pr}(X\le x)\le1/2\\
-\log_{2}(2\mbox{Pr}(X\ge x)) & \mbox{if }\mbox{Pr}(X\le x)>1/2
\end{array}.\right.
\]
We call this the $\zeta$-score because, like the $z$-score, it measures
the extent to which an observation is extreme, or in the tails of
the distribution. For the $z$-score the unit is the standard deviation
while for the $\zeta$-score it is the tail half-area: for $0<\zeta_{1}<\zeta_{2}$
such that $\zeta_{2}-\zeta_{1}=1,$ the tail area for $x_{2}$ is
half that for $x_{1}$. For normal data the $z$-score can be converted
to tail areas while the $\zeta$-score describes directly the tail
areas for any distribution: $-1,0$, and $1$ correspond to Q1, median,
and Q3, respectively. Since $|\zeta|=4$ corresponds to $2T\!A=.0625$,
$\zeta$-scores from 3 to 6 are portions of the tail that generally
are of inferential interest. The $\zeta$-scores for each of the estimators
were calculated and plotted against those obtained from the mean so
that the line with slope 1 represents the sample mean. The curve with
the smallest slope is the median and the curve with slope between
the mean and median is the $t_{3}$ mle. 

The information efficiencies are 0.724 for the median and 0.906 for
the $t_{3}$ mle. Since the sample size is $n=10$ this corresponds
to a loss of a little less than 3 observations and about 1 observation,
respectively. The dotted lines are for the sample mean based on samples
of size $n=7$ and $n=9$. 

For the graph on the right, the above was repeated using $t_{3}$
data rather than normal data. The sample mean is still used as a reference
so the line with slope 1 is for the sample mean. The sample mean is
less efficient than the other estimators for $t_{3}$ data so that
their curves have slope greater than 1. 

Using the variances (obtained from 100,000 samples), the relative
efficiency compared to the sample mean is 1.50 and 1.78 for the median
and $t_{3}$ mle, respectively. The dotted lines are for the sample
mean based on samples of size $n=15$ and $n=17$. Clearly, interpretation
of efficiency in terms of sample size does not work here. Caution
should be exercised when data do not follow the normal distribution. 

Instead of changing the sample sizes, the comparison of estimators
in terms of slope suggests rescaling the data. The dashed lines show
the sample mean when the $t_{3}$ data are multiplied by $1.50^{-1/2}$
and $1.78^{-1/2}$. These dashed lines are much closer to the median
and $t_{3}$ mle than those based on sample size adjustments. For
normal data, the result of rescaling is identical to changing the
sample size. For example, the distribution of the sample mean when
$n=9$ is the same when each observation is scaled by $\sqrt{10/9}$
and $n=10$. 

\begin{figure}
\centering

\includegraphics[width=0.5\columnwidth]{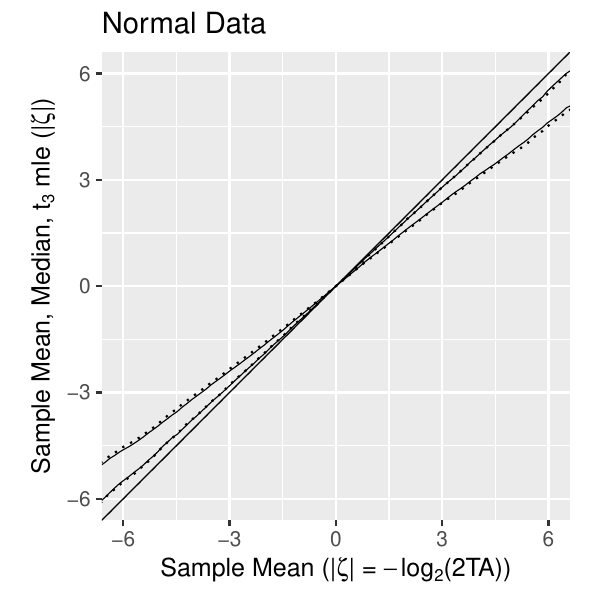}\includegraphics[width=0.5\columnwidth]{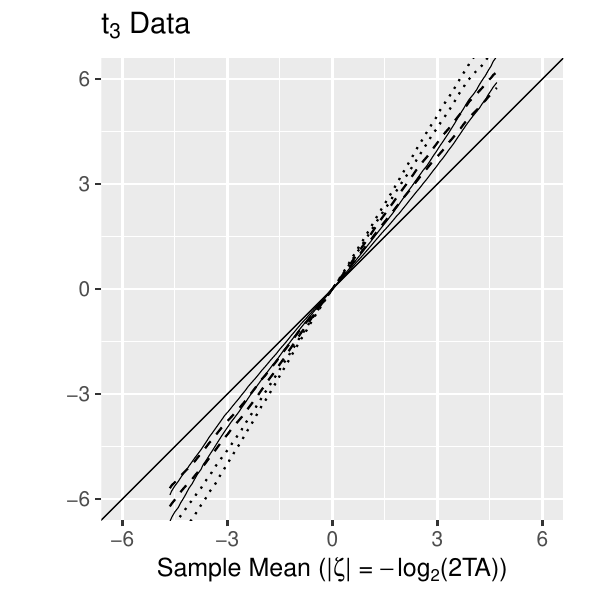}

\caption{\label{fig:Comparing-Estimators:-Mean,}Comparing Estimators: Mean,
Median, $t_{3}$ mle.}
\end{figure}

\subsection{Inference for an Odds Ratio\label{subsec:ExampleInference-for-an}}

Consider an inference problem comparing two populations with dichotomous
data and independent samples. Assume the sample space is $\mathcal{Y}=\{0,\cdots,n_{1}\}\times\{0,\cdots,n_{2}\}$,
where $n_{1},n_{2}\ge1$ are the two sample sizes, and the data are
generated from two Binomial distributions $x_{1}\sim Bin(n_{1},p_{1})$
and $x_{2}\sim Bin(n_{2},p_{2})$, where $p_{1},p_{2}\in(0,1)$ are
the two population success probabilities. We define the parameter
of interest to be the log odds ratio 
\[
\theta=\log\left(\frac{p_{1}}{1-p_{1}}\right)-\log\left(\frac{p_{2}}{1-p_{2}}\right)\in\mathbb{R}
\]
between the two populations and let $\undertilde{\theta}=n_{1}p_{1}+n_{2}p_{2}\in(0,n_{1}+n_{2})$
be the nuisance parameter. We illustrate the benefits of using this
nuisance parameterization versus other choices in the next, although
nuisance parameterization should not affect the inference on the parameter
of interest $\theta$.

The log-likelihood function of the data $y=(x_{1},x_{2})$ in the
success probabilities $(p_{1},p_{2})$ can be written as 
\begin{equation}
\ell_{y}(p_{1},p_{2})\propto x_{1}\log p_{1}+(n_{1}-x_{1})\log(1-p_{1})+x_{2}\log p_{2}+(n_{2}-x_{2})\log(1-p_{2}).\label{eq:logl}
\end{equation}
Taking derivatives of the log-likelihood (\ref{eq:logl}) with respect
to $\theta$ and $\undertilde{\theta}$, we get the two score functions
\begin{eqnarray}
s_{y}(\theta,\undertilde{\theta}) & = & \frac{x_{1}-n_{1}p_{1}}{p_{1}(1-p_{1})}\nabla p_{1}+\frac{x_{2}-n_{2}p_{2}}{p_{2}(1-p_{2})}\nabla p_{2}\label{eq:score1}\\
\undertilde{s}_{y}(\theta,\undertilde{\theta}) & = & \frac{x_{1}-n_{1}p_{1}}{p_{1}(1-p_{1})}\widetilde{\nabla}p_{1}+\frac{x_{2}-n_{2}p_{2}}{p_{2}(1-p_{2})}\widetilde{\nabla}p_{2}\label{eq:score2}
\end{eqnarray}
where $s_{y}(\theta,\undertilde{\theta})=s(y,\underline{\theta})$
and $\undertilde{s}_{y}(\theta,\undertilde{\theta})=\undertilde{s}(y,\text{\ensuremath{\underbar{\ensuremath{\theta}}}})$.
The covariance of these two score vectors is 
\begin{eqnarray*}
\mathsf{E}\left[s\undertilde{s}\right] & = & \frac{n_{1}\nabla p_{1}\widetilde{\nabla}p_{1}}{p_{1}(1-p_{1})}+\frac{n_{2}\nabla p_{2}\widetilde{\nabla}p_{2}}{p_{2}(1-p_{2})}.
\end{eqnarray*}

However, because of
\begin{equation}
\frac{\widetilde{\nabla}p_{1}}{p_{1}(1-p_{1})}-\frac{\widetilde{\nabla}p_{2}}{p_{2}(1-p_{2})}=0\label{eq:delta1}
\end{equation}
and
\begin{equation}
n_{1}\nabla p_{1}+n_{2}\nabla p_{2}=0,\label{eq:delta2}
\end{equation}
it is easy to show that $E[s\undertilde{s}]=0$ and 
\[
\undertilde{s}=\frac{(x_{1}+x_{2}-\undertilde{\theta})\widetilde{\nabla}p_{1}}{p_{1}(1-p_{1})}.
\]
 Equation (\ref{eq:delta1}) is the result of the interest parameterization,
while (\ref{eq:delta2}) is the result of the nuisance parameterization.
The two benefits of this nuisance parameterization are that
\begin{enumerate}
\item the two scores are orthogonal to each other, i.e., $\mathsf{E}[s\undertilde{s}]=0$,
and
\item the solution of the nuisance parameter $\hat{\undertilde{\theta}}=x_{1}+x_{2}$
to $\undertilde{s}=0$ is independent of the interest parameter $\theta$.
\end{enumerate}
The first benefit means that there is no need to orthogonalize the
score $s$, while the second benefit makes the inference using profiled
score easier.

Because of (\ref{eq:delta2}), the standardized score of $s$ can
be written as
\begin{equation}
\bar{s}_{y}(\theta,\undertilde{\theta})=\left(\frac{1}{n_{1}p_{1}(1-p_{1})}+\frac{1}{n_{2}p_{2}(1-p_{2})}\right)^{-\frac{1}{2}}\left(\frac{\bar{x}_{1}-p_{1}}{p_{1}(1-p_{1})}+\frac{\bar{x}_{2}-p_{2}}{p_{2}(1-p_{2})}\right)\label{eq:stdscore}
\end{equation}
with information
\[
\mathsf{E}\left[\nabla\bar{s}\right]^{2}=\left(\frac{1}{n_{1}p_{1}(1-p_{1})}+\frac{1}{n_{2}p_{2}(1-p_{2})}\right)^{-1}.
\]
The information is obtained by noting
\[
\frac{\nabla p_{1}}{p_{1}(1-p_{1})}-\frac{\nabla p_{2}}{p_{2}(1-p_{2)}}=1.
\]

For inference of the interest parameter, the $z$-standard intervals
can be defined as
\[
CI_{-z}(y)=\{\theta:\bar{s}_{y}(\theta,\undertilde{\theta})\ge-z\},\:CI_{+z}(y)=\{\theta:\bar{s}_{y}(\theta,\undertilde{\theta})\le z\},
\]
and the two-sided interval can be obtained by the intersection $CI_{z}(y)=CI_{-z}(y)\cap CI_{+z}(y)$.
Although the choice of nuisance parameterization does not affect the
standardized score $\bar{s}$ and the interval estimation, they depend
on the value of the nuisance parameter $\undertilde{\theta}$. A popular
method is to use the profiled score $\bar{s}_{y}(\theta,\undertilde{\hat{\theta}})$
by plugging in the profiled estimate of the nuisance parameter, i.e.,
$\hat{\undertilde{\theta}}=x_{1}+x_{2}$ in this case. Having $\hat{\undertilde{\theta}}$
independent of the interest parameter $\theta$ simplifies the process
of finding and assessing the $z$-standard intervals. In finding the
$z$-standard intervals, it is easier to solve $\bar{s}_{y}(\theta,\undertilde{\theta})\ge-z$,
$\bar{s}_{y}(\theta,\undertilde{\theta})\le z$, or $\bar{s}_{y}(\theta,\undertilde{\theta})^{2}\le z^{2}$
for $p_{1},p_{2}\in(0,1)$ by conditioning on $n_{1}p_{1}+n_{2}p_{2}=x_{1}+x_{2}$.
Then, for each pair of estimates $\hat{p}_{1}$ and $\hat{p}_{2}$,
we can calculate the odds ratio $\widehat{or}=(\hat{p}_{1}(1-\hat{p}_{2}))/((1-\hat{p}_{1})\hat{p}_{2})$
or equivalently the log odds ratio $\log(\widehat{or})$. Boundary
issues when $x_{1}=0$ or $n_{1}$ or when $x_{2}=0$ or $n_{2}$
may need some special attention, and whether to include the equal
sign in the definitions of the $z$-standard intervals makes a difference
especially when the true success probabilities are close to 0 or 1.

Before we give the algorithm to find the $z$-standard intervals,
we are interested to know the coverage of the $z$-standard intervals
and its relation to the choice of the nuisance parameter value. Although
profiling suggests $\hat{\undertilde{\theta}}=x_{1}+x_{2}$ as the
nuisance parameter value, other values can be attempted. For example,
an adjustment to the profiled nuisance parameter value can be made
by 
\[
\hat{\undertilde{\theta}}(c)=\frac{n_{1}(x_{1}+c)}{n_{1}+2c}+\frac{n_{2}(x_{2}+c)}{n_{2}+2c},
\]
where $c>0$ is a positive constant, to mimic the plus-4 estimate
of a proportion. This conservative adjustment can help with the boundary
cases in computation by avoiding parameter values outside the parameter
space. We conducted numerical studies with $n_{1}$ = 20 and $n_{2}$
= 30 for various odds ratio setups and adjustments to the nuisance
parameter value in Table 1. It can be seen that 1) the $z$-stardard
intervals have coverage close to the nominal level when the true success
probabilities are not too close to 0 or 1, 2) the coverage is a non-decreasing
function of the adjustment constant $c$, 3) there is a difference
between the options with and without the equal sign in the definition
of the $z$-standard interval when the success probabilities are close
to 0 or 1, but this difference can be mitigated with an adjustment.
In comparison, Table 1 also include the coverage of Fisher's exact
intervals. 
\begin{table}
\centering

\begin{tabular}{|c|c|c|c|c|c|c|c|c|c|c|}
\hline 
\multicolumn{3}{|c|}{True parameters} & \multicolumn{4}{c|}{With ``=''} & \multicolumn{4}{c|}{Without ``=''}\tabularnewline
\hline 
OR & p1 & p2 & c = 0 & c = 0.5 & c = 1 & Fisher & c = 0 & c = 0.5 & c = 1 & Fisher\tabularnewline
\hline 
\hline 
1.0 & 0.01 & 0.01 & 0.999 & 1.000 & 1.000 & 1.000 & 0.394 & 1.000 & 1.000 & 1.000\tabularnewline
\hline 
1.0 & 0.20 & 0.20 & 0.957 & 0.958 & 0.960 & 0.980 & 0.957 & 0.958 & 0.960 & 0.980\tabularnewline
\hline 
1.0 & 0.50 & 0.50 & 0.944 & 0.944 & 0.944 & 0.975 & 0.944 & 0.944 & 0.944 & 0.975\tabularnewline
\hline 
1.0 & 0.70 & 0.70 & 0.954 & 0.954 & 0.954 & 0.977 & 0.954 & 0.954 & 0.954 & 0.977\tabularnewline
\hline 
1.0 & 0.90 & 0.90 & 0.968 & 0.976 & 0.981 & 0.989 & 0.962 & 0.976 & 0.981 & 0.989\tabularnewline
\hline 
1.5 & 0.015 & 0.01 & 1.000 & 1.000 & 1.000 & 1.000 & 0.452 & 1.000 & 1.000 & 1.000\tabularnewline
\hline 
1.5 & 0.273 & 0.20 & 0.946 & 0.958 & 0.958 & 0.981 & 0.946 & 0.958 & 0.958 & 0.981\tabularnewline
\hline 
1.5 & 0.60 & 0.50 & 0.946 & 0.946 & 0.946 & 0.976 & 0.946 & 0.946 & 0.946 & 0.976\tabularnewline
\hline 
1.5 & 0.778 & 0.70 & 0.952 & 0.953 & 0.953 & 0.975 & 0.952 & 0.953 & 0.953 & 0.975\tabularnewline
\hline 
1.5 & 0.931 & 0.90 & 0.959 & 0.982 & 0.982 & 0.995 & 0.959 & 0.982 & 0.982 & 0.995\tabularnewline
\hline 
4.0 & 0.039 & 0.01 & 0.982 & 0.998 & 0.998 & 0.998 & 0.647 & 0.998 & 0.998 & 0.998\tabularnewline
\hline 
4.0 & 0.50 & 0.20 & 0.952 & 0.957 & 0.957 & 0.977 & 0.952 & 0.957 & 0.957 & 0.977\tabularnewline
\hline 
4.0 & 0.80 & 0.50 & 0.948 & 0.948 & 0.950 & 0.978 & 0.948 & 0.948 & 0.950 & 0.978\tabularnewline
\hline 
4.0 & 0.903 & 0.70 & 0.955 & 0.960 & 0.970 & 0.987 & 0.955 & 0.960 & 0.970 & 0.987\tabularnewline
\hline 
4.0 & 0.973 & 0.90 & 0.938 & 0.974 & 0.974 & 0.989 & 0.914 & 0.974 & 0.974 & 0.989\tabularnewline
\hline 
\end{tabular}

\caption{Coverage of $z$-standard interval vs. Fisher's exact internal for
sample sizes $n_{1}=20$ and $n_{2}=30$. The nominal coverage is
set at 0.95.}
\end{table}

\begin{figure}
\centering

\includegraphics[height=0.22\paperheight]{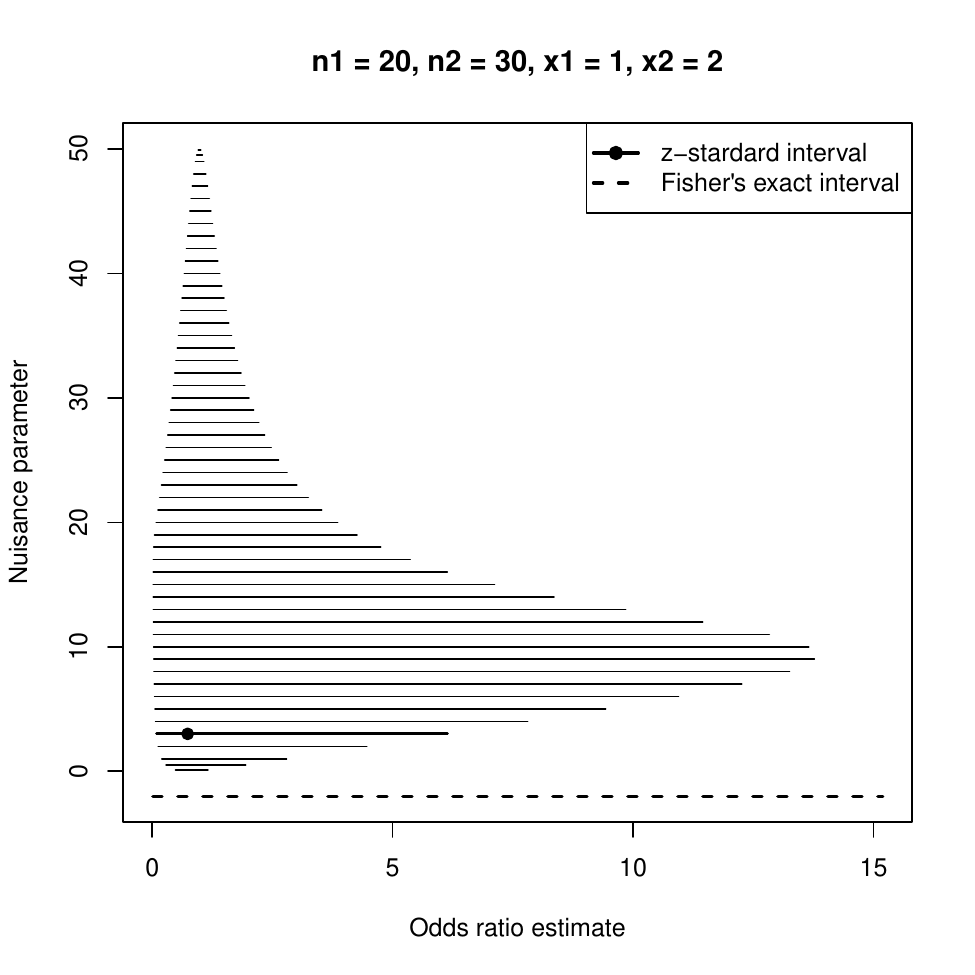}

\includegraphics[height=0.22\paperheight]{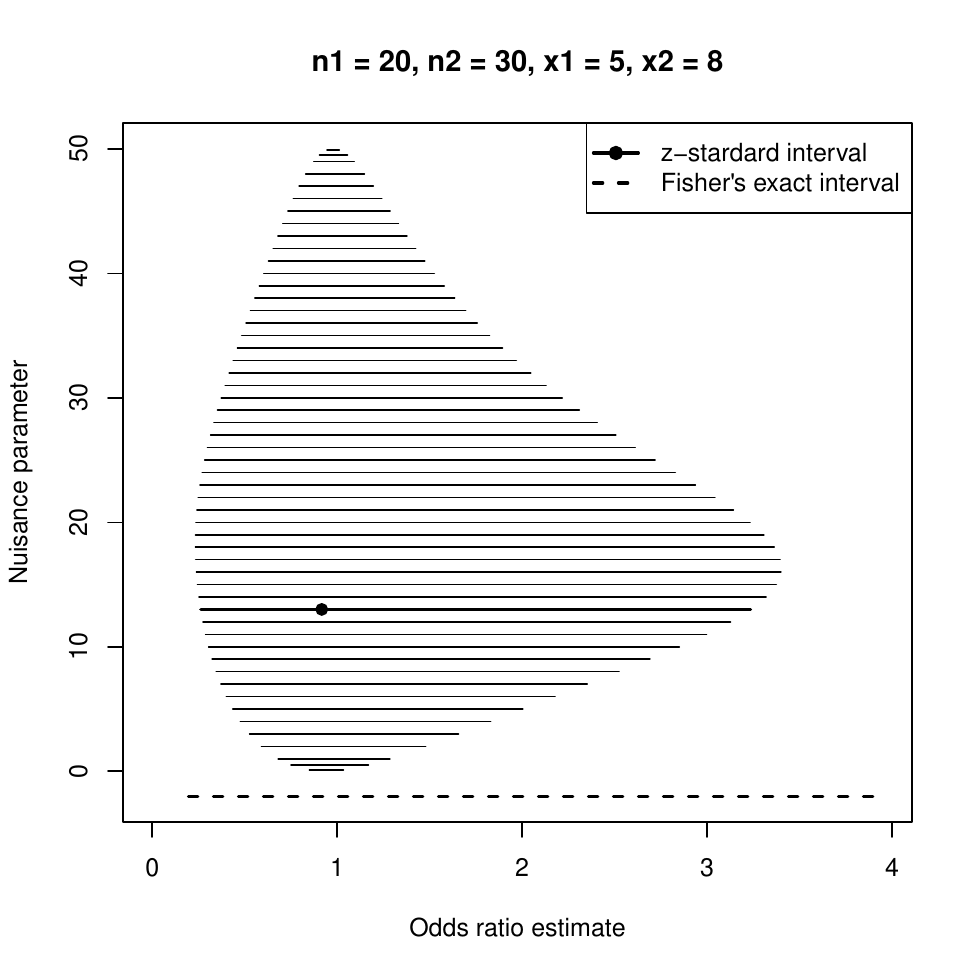}

\includegraphics[height=0.25\paperheight]{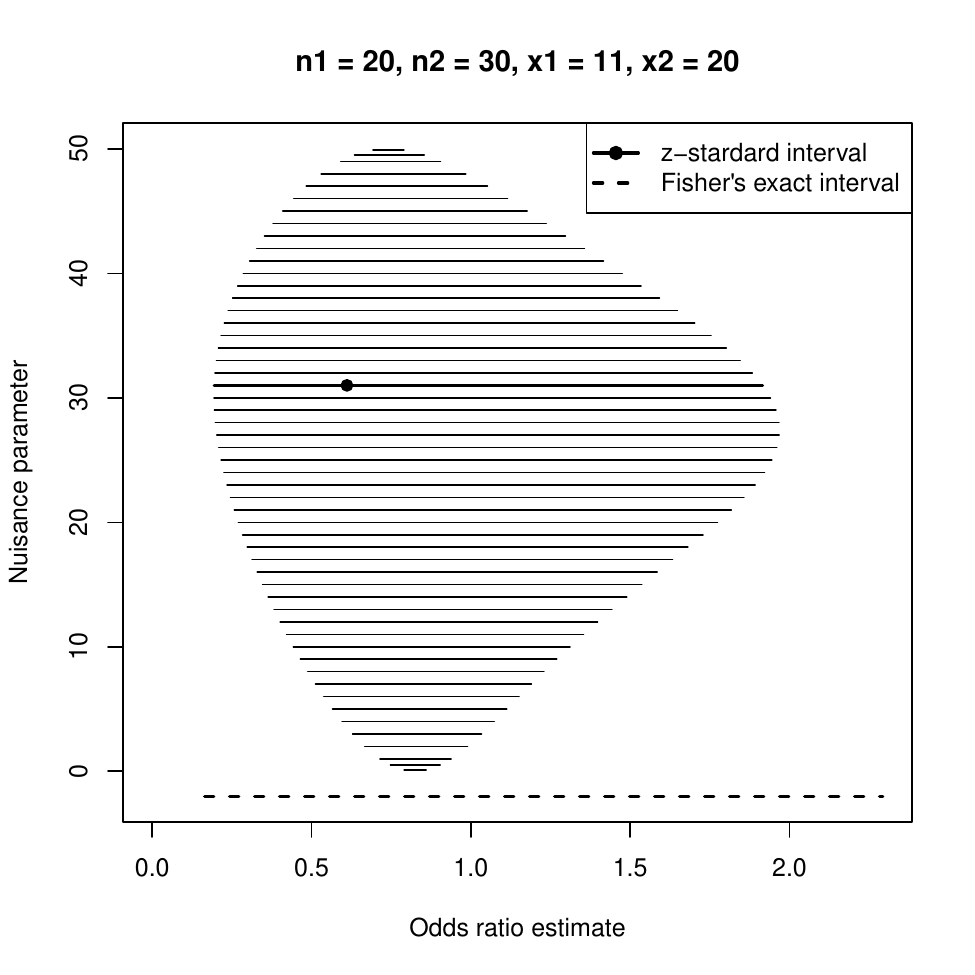}

\caption{\label{fig:Three-examples-of}Three examples of $z$-standard interval
(at various nuisance parameter values) vs. Fisher's exact interval.}
\end{figure}

Next, we give an algorithm to find the two-sided $z$-standard intervals.
The one-sided intervals can be found in a similar way. The two-sided
interval can be found by solving $\bar{s}_{y}(\theta,\undertilde{\theta})^{2}\le z^{2}$,
which can be written as a six-degree polynomial inequality in $p_{1}$
given $\undertilde{\theta}=n_{1}p_{1}+n_{2}p_{2}$. The polynomial
is facing upward, so the six solutions give up to three intervals,
among which the middle one is the desired solution given the constraints.
Figure \ref{fig:Three-examples-of} shows three examples with $z$-standard
intervals for the odd ratio with $\undertilde{\hat{\theta}}=x_{1}+x_{2}$
compared to those at various nuisance parameter values, as well as
the Fisher's exact interval. 

\begin{figure}
\centering\includegraphics[width=0.7\columnwidth]{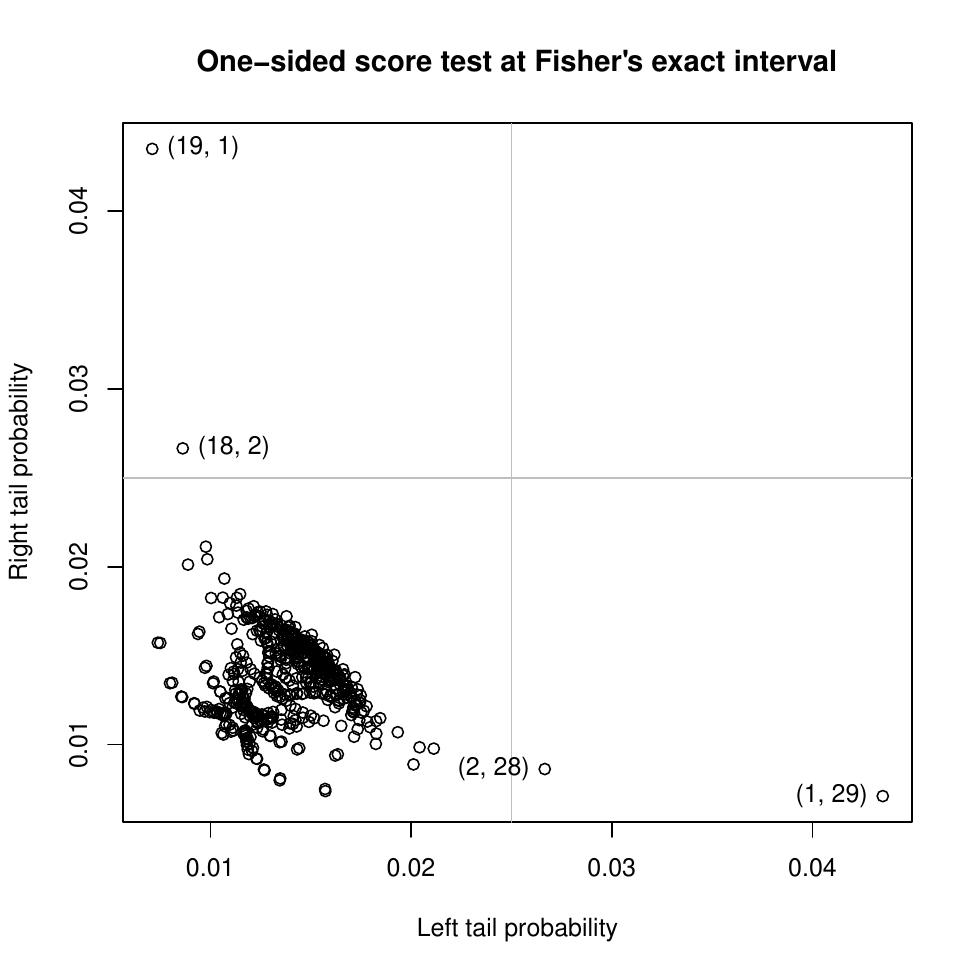}

\caption{\label{fig:Tail-probabilities-of}Tail probabilities of one-sided
score test at the endpoints of the Fisher's exact interval with a
0.95 nominal level. The tail probabilities are given for $n_{1}=20$,
$n_{2}=30$, $x_{1}=1,\cdots,19$, and $x_{2}=1,\cdots,29$.}
\end{figure}

Table 1 shows that the coverage of the Fisher's exact interval is
higher than that of the $z$-standard interval and above the nominal
level. Thus, it would be interesting to know if the Fisher's exact
interval is wider than the $z$-standard interval as suggested in
Figure \ref{fig:Three-examples-of}. However, instead of directly
comparing the width of the intervals, we find the tail probability
of one-sided score test using (\ref{eq:stdscore}) at the endpoints
of the Fisher's exact interval. The results are showing in Figure
\ref{fig:Tail-probabilities-of}. As can be seen, both the left and
right tail probabilities are under the nominal 0.025 level, except
for 4 extreme cases. So we can conclude that the Fisher's exact interval
is in general wider than the $z$-standard interval and maybe more
shifted to one side in extreme cases. 

\section{Discussion\label{sec:Discussion}}

Generalized estimation is related to, but distinct from, the theory
of estimating functions. A generalized estimate $g$ can serve as
an estimating function in that its root provides a point estimate.
As noted above, requirements on the function $g$ are slightly different
from those in \citet{Godambe1960} and the results are stated for
information rather than variance. An important distinction is that
in generalized estimation $g$ \emph{is} the estimator rather than
a vehicle for finding an estimator. This distinction is illustrated
in the examples where we use $g$ directly for calculating a confidence
interval or assessing how far an observation is in the tails of the
distribution ($\zeta$-value). For fixed $y\in\mathcal{Y}$, the estimate
$g_{y}$ indicates the model most consistent with $y$ by assigning
the value zero and assigns real values to other models to indicate
their consistency, or lack thereof, with $y$. When $\theta_{\circ}\in\Theta$
is fixed, $g$ is a test statistic for $H_{\circ}:\theta=\theta_{\circ}$.
The difference from typical hypothesis testing is that we consider
a continuum of hypothesis and estimators are assessed by how they
change with $\theta_{\circ}$ and so their ability to use the data
to distinguish models in $\Theta$. 

\citet{LiangZeger1995} indicate two limitations of estimating functions. 
\begin{quote}
First, optimality is ascribed to the estimating function, but scientists
and other practitioners are concerned about estimators. As Crowder
(1989) put it: ``This is like admiring the pram rather than the baby.''
\end{quote}
Given the inferential value of the function $g$ beyond simply reporting
its root (point estimate) a more apt description would be that point
estimation is like admiring the baby's rattle rather than the baby.
The second limitation concerns nuisance parameters and we recognize
this as a difficulty for many forms of inference. Requiring generalized
estimators to be orthogonal to the nuisance space address some of
these problems. 

Another difference is that estimating functions are used to provide
new point estimates for both parametric and semi-parametric inference.
A major result of this paper is not the introduction of new estimators
but the superiority of inference based on the score and log likelihood
function. The information $\Lambda$ provides a stronger justification
for their superiority than do assessments of point estimators such
as variance and bias. Of the extensive estimating functions literature,
\citet{Chan1998} is the closest to our work. Key differences are
that these authors do not consider the connection of the parameter
space with the space of probability distributions, and so do not address
the parameter invariance of $g$ and their approach to nuisance parameters
will not be invariant under reparameterization. 

The inference problem in this paper is that of describing what a sample
$y$ tells us about the possible models in a family of distributions.
The plausibility of each model is assessed by calculating how extreme
$y$ would be for the given model. The notion of 'extreme' is made
precise by the generalized estimator $g$ which for each model provides
an equivalence relation on $\mathcal{Y}$ and a total ordering of
these equivalence classes. For $y_{1},y_{2}\in\mathcal{Y}$, $y_{1}\sim_{\underline{\theta}}y_{2}$
if-f $g(y_{1},\underline{\theta})=g(y_{2},\underline{\theta})$ and
$\left[y_{1}\right]\preccurlyeq_{\underline{\theta}}\left[y_{2}\right]$
if-f $g(y_{1},\underline{\theta})\le g(y_{2},\underline{\theta})$
where $\left[y\right]=\left[y\right]_{\underline{\theta}}$ is the
equivalence class for $y$. The numerical value associated with this
notion of extreme sample is the tail area, or on the log scale, the
$\zeta$-score. For the tail area to be considered a probability requires
that the study data were obtained according to the sampling plan used
to define the models in $M$. In many cases, this means $x_{1}^{obs},x_{2}^{obs},\ldots,x_{n}^{obs}$
were obtained by a simple random sample. There is no need to consider
hypothetical repeated samples. \citet{Vos2022} provide further discussion
on the role of hypothetical repeated trials in frequentist inference. 

A related, but distinct, inference problem is that of using $y$ to
obtain a point estimate that is close to the 'true' distribution.
Here 'true' means the member in the family of distributions that is
the best approximation to the distribution of the population of interest.
As the true distribution is unknown, inference procedures are described
for a generic distribution in the family. This focus on a single distribution
is the key distinction from our approach. The tools for assessing
point estimators, mean, variance, and mean square error, are all defined
pointwise. As the distribution considered is generic, these are defined
on the entire parameter space, but each depends only on the distribution,
not on models in the neighborhood of the distribution. In contrast,
Fisher information describes the relationship between models in the
family, heuristically the amount of information in the sample to distinguish
between nearby models. Theorems \ref{thm:1} and \ref{thm:Lambda}
can be described as versions of the Cram\textipa{\'{e}}r-Rao bound
that recognize this distinction between information and variance. 

The standard normal distribution, $m_{\circ}=N(0,1)$, provides a
simple example of this distinction. The sample mean $\bar{x}$ has
expectation 0 and variance $n^{-1}$ at $m_{\circ}$.The Fisher information
is not defined for $m_{\circ}$. Fisher information is defined for
a smooth family of distributions. We consider two families that contain
$m_{\circ}$. For $M_{\mathcal{X}_{1}}=\left\{ m:m=N(\mu,1),\mu\in\mathbb{R}\right\} $
the Fisher information for $y$ is $n$ while for $M_{\mathcal{X}_{2}}=\left\{ m:m=N(0,\sigma),\ 0<\sigma\right\} $
it is $2\sigma^{-2}n$ or $2n$ at $m_{\circ}$. Fisher information
describes the family, not the estimator. For $M_{\mathcal{X}_{1}}$
the information utilized by the sample mean is $\Lambda_{1}(\bar{x})=n$
while for $M_{\mathcal{X}_{2}}$ it is $\Lambda_{2}(\bar{x})=0$ corresponding
to the fact that a location estimator is of no inferential value when
all the distributions in the family have the same mean. 

Some authors require that a parameter have a meaningful interpretation
\citep{McCullagh2002}. The role of parameters in our approach is
to describe the smooth structure of $M$ and estimators on $M$, and
for this, the only required property is that the parameterization
be a diffeomorphism. When deciding whether a family of distributions
$M$ are appropriate for real-world data, it is useful to consider
the attributes of specific distributions (points in $M$) in relation
to the sample space $\mathcal{X}$. Common attributes are the mean
and variance. Often, but not necessarily, the attribute function is
smooth and so a parameterization in the sense used in this paper.
An important distinction between a parameterization and an attribute
function is that the latter often requires units for it to be meaningful.
For example, if the data are measured in grams, then, as an attribute,
the mean is also expressed in grams. For generalized estimation, if
the mean is used as a parameterization these units are ignored because
a generalized estimate does not take values in the parameter/attribute
space but is a function on this space. The graph of $g$ may have
units on the abscissa but there are no units on the ordinate. Units
are not required for the equivalence relation and ordering described
above. 

Point estimation can be done without using any parameterization by
replacing the estimate $\hat{\theta}$ with $\hat{m}$, the distribution
named by $\hat{\theta}$, and replacing Euclidean distance on $\underline{\Theta}$
with Kullback Leibler divergence on $M$ to define distribution versions
of the mean, variance, and mean square error. However, these quantities,
like their parametric versions, are defined point-wise and so do not
capture the structure of the family of distributions. Details can
be found in \citet{Wu2012} and \citet{Vos2015}.

We have addressed nuisance parameters by orthogonalizing the estimator.
Another approach is to use an orthogonal parameterization. See, for
example, \citet{CoxReid1987}. As the example in Section \ref{subsec:ExampleInference-for-an}
shows, these can be convenient for calculations. We chose to orthogonalize
estimators because orthogonal parameters need not exist (unless the
interest parameter is a scalar) and we want to de-emphasize parameters
as we are making inferences for distributions rather than select attributes
of distributions. The orthogonalized Fisher information for the interest
parameter has been defined previously by \citet{Liang1983-ve} and
\citealt[(][ pp. 250-252)]{Amari1990-bo}. \citet{Hudson2000} describe
the role of conditional and marginal orthogonalized Fisher information
for inference in presence of nuisance parameters. 

The $\zeta$-score used in Section \ref{subsec:ExampleNormal-and--distributions}
is proportional to the transformation used in Fisher's method to combine
p-values. Twice the tail area, used in the definition of $\zeta$,
is the p-value for the test conditional on the observed tail where
each tail has probability $1/2$.

\section*{}

\bibliographystyle{apalike2}
\bibliography{../vos}

\end{document}